\documentstyle{article}

\newcommand{\N}{{\bf N}}
\newcommand{\Z}{{\bf Z}}
\newcommand{\bP}{{\bf P}}
\begin{document}

\title{Vanishing theorems for products of exterior and
symmetric powers}
\author{F. Laytimi\thanks{\it Mathematics Department, University of
Lille}
\and W. Nahm\thanks{\it Physikalisches Institut, Bonn University}}

\maketitle
\section{Introduction}
\setcounter{page}{1}

Many problems of algebraic geometry involve
vanishing theorems for  cohomology groups of ample
vector bundles.

For ample vector bundles $E$ over compact complex
varieties $X$ and a Schur functor $S_I$ corresponding
to an arbitrary  partition $I$ of the integer $|I|$,
one would like to know the optimal vanishing theorem for the
cohomology groups $H^{p,q}(X, S_I(E))$, depending on the rank of $E$
and the dimension $n$ of $X$. For $p=n$, a conjectural
vanishing statement has been motivated by the study
of the connectedness of degeneracy loci [La].

In an unpublished paper one of us
proved a vanishing theorem for the situation where the partition
$I$ is a hook [Na]. Here we give a
simpler proof of this theorem. In Theorem 2.2 we treat the same
problem under weaker positivity assumptions, in particular
under the
hypothesis of ample $\Lambda ^m E$  with $m\in \N^*$.
In this case we also need some bound on the weight $|I|$ of the
partition.
Moreover, we prove that the same vanishing condition
applies for $H^{q,p}(X, S_I(E))$, with $p,q$ interchanged.

\section{Statement and Results}

Let $E$ a vector bundle of rank $e$.
For $ 0\leq \alpha< k$ the hook Schur functor ${\Gamma}^\alpha_k E$,
corresponding to the partition
$(\alpha  +1,1,\dots ,1)$ of $k\in \N^*$,
can be defined inductively as follows:
$$\Gamma^0_k E=\Lambda^kE$$
and
$$\Lambda^{k-\alpha }E\otimes S^{\alpha } E=
\Gamma^{\alpha-1}_k E\oplus \Gamma ^\alpha_k E$$
for $0<\alpha<k$.
In particular,  $\Gamma^{k-1}_k E=S^kE$.
Note that $\Gamma ^\alpha_k E=0$\\ for $e-k+\alpha<0$.

Define a function $\delta:\N\rightarrow  \N^*$ by:
$${\delta(x)\choose 2}\le x<{\delta(x)+1\choose 2}\ .$$
In particular, $\delta({m\choose 2})=m$ and $\delta(0)=1$.
\medskip

{\bf Theorem 2.1 }(Nahm 1995) {\bf :} {\it Let $X$ be a compact 
complex variety of dimension $n$
and $E$ an ample vector bundle over $X$ of rank $e$.
 Then,
$$H^{p,q}(X,\Gamma ^\alpha_k E)=0,$$ for
$q+p-n>(\delta(n-p)+\alpha)(e-k+2\alpha)- \alpha(\alpha+1)\ .$}
\bigskip

For $\alpha=k-1$ one obtains
$H^{p,q}(X,S^\alpha E\otimes det E)=0,$
when
$q+n-p>(\delta (n-p)-1)k$.
For $p=n$, this specialises to Griffith's vanishing theorem [Gr].

For $\alpha=0$ one obtains
$H^{p,q}(X,\Lambda^k E)=0$ for
$q+p-n>\delta(n-p)(e-k)\ .$
For $p=n$, this specialises to Le Potier's vanishing theorem [LeP].
\bigskip

Here we prove the more general
\smallskip

{\bf Theorem 2.2 :} {\it  Let $X$ be a compact complex variety of
dimension $n$.
Let $k=ml+s$ with $k,m,l\in \N^*$, $0\leq s<m$.
Let $E$ be a vector bundle over $X$ of rank $e$,  such that
$S_{I(l,m,s)}E$ is ample, with the partition of $k$
$$I(l,m,s)=(\ \underbrace {l+1,\ldots ,l+1}_{s\mit\ times}\ ,
\underbrace {\ l, \ldots ,l\ }_{(m-s)\mit\ times})\ .$$
Let $ 0\leq \alpha<k$,\  $0<p\leq n$, and
$$r=\delta \left(n-p+{m\choose 2}\right)\ .$$
Suppose that one of the following four conditions is true:

-- $m=1$

-- $\alpha=0$, \  $(r-1)k>r(n-p-1)$

-- $\alpha\geq 1$, \ $(r+\alpha-2)k>
(r+\alpha-1)\left(n-p+ra+{\alpha\choose 2}-1\right)+
\delta_{\alpha,2}\delta_{m,2}\delta_{n-p,1}$

\qquad (with Kronecker's $\delta$-notation)

-- $k>e+1$, $(r+\alpha-\beta-1)k>(r+\alpha-\beta)
\left(n-p+ra+{\alpha\choose 2}-{\beta\choose 2}-1\right)$,

\qquad where $\beta=k-e$.

Then

$$H^{p,q}(X,\Gamma ^\alpha_k E)=H^{q,p}(X,\Gamma ^\alpha_k E)=0,$$ 
when
$$q+p-n> (r +\alpha)(e-k+\alpha)+\alpha(r-1)\ .$$}\\
\smallskip
{\bf Remark 2.3 :}
For $m=1$, the condition $S^lE$ ample is equivalent to $E$ ample,
such that we obtain theorem 2.1 with the additional freedom to
interchange $p$ and $q$.
When $k$ is divisible by $m$, the ampleness condition
is equivalent to $\Lambda^m E$ ample, as we shall show now.
At the same time, we shall see that theorem 2.2 only yields
vanishing statements for ample vector bundles, as expected.
\medskip

Recall the dominance partial order of
partitions\\
 $I=(i_1,i_2\dots )$ of $|I|$ and $J=(j_1,j_2\dots )$ of $|J|$.
Usually it is defined under the
assumption that $I,J$ are partitions of the same integer, but we
generalise it to arbitrary pairs of partitions by scaling 
$I$ and $J$ to rational partitions
of 1. Let
$\stackrel{-}{i}_k=i_k/{\vert I\vert}$,
\hskip2mm $\stackrel{-}{j}_k=j_k/{\vert J\vert}$.
\smallskip
We say that $S_I$ dominates $S_J$, if
$\Sigma _{k=1}^l \stackrel{-}{i}_k\geq \Sigma _{k=1}
^l \stackrel{-}{j}_k$ for any $l$.
If each one dominates the other, i.e. if $I$ and $J$ are
proportional, we call these Schur functors
equivalent. For example, $S_{I(l,m,s)}$ dominates 
$\Gamma ^\alpha_k$ and $S_{I(l,m,0)}$ is equivalent to 
$\Lambda^m$.
\bigskip

{\bf Lemma 2.4 :} {\it If $S_I$ dominates $S_J$, then ampleness 
of $S_IE$ implies ampleness of $S_J E$.}
\bigskip

{\it Proof:}
It is sufficient to show that every irreducible direct summand of
$S_J E^{\otimes k{\vert I\vert}}$
is isomorphic to a direct summand
of $S_I E^{\otimes k{\vert J\vert}}$
for some $k>0$. More generally, consider the set $T'$ of
partitions $I'$ with $|I'|={\vert J\vert}k{\vert I\vert}$
which are dominated by $k{\vert J\vert}I$. It is sufficient to show
that for any $I'\in T'$
the bundle $S_{I'}E$ is isomorphic to a direct summand
of $S_I E^{\otimes k{\vert J\vert}}$
when $k$ is a (sufficiently large) multiple of $lcm(1,2,\ldots,e)$.
In this case, all vertices of the convex hull of $T'$ in 
${\bf R}^e$
different from $({\vert J\vert}k{\vert I\vert}/e,\ldots,
{\vert J\vert}k{\vert I\vert}/e)$ lie on a face $F_l$ passing through
$k{\vert J\vert}I$. Such a face is given by the equation
$\Sigma _{n=1}^l \stackrel{-}{i}_n= \Sigma _{n=1}
^l \stackrel{-}{i'}_n$ for some $l<e$.

For partitions $I$ of length at most $e$
let $V(I)=S_I {\bf C}^e$ be the corresponding highest weight
representation of $GL_e({\bf C})$. Define
$$T_{e,r}=\{(I_1,\ldots,I_r,I')\
\vert\  V(I')\subset V(I_1)\otimes \ldots
\otimes V(I_r)\}\ .$$
Recently, Knutson and Tao
have shown that the semigroups $T_{e,r}$ are saturated, such that
for fixed $I_1,\ldots,I_r$ the set
$\{I'\vert (I_1,\ldots,I_r,I')\in T_{e,r}\}$ is convex. More
precisely, it contains all integral points of its convex hull $H$ in
${\bf R}^e$ [K-T] (see also [Bu] and [Ze]).

Let us specialize to $r=k{\vert J\vert}$ and $I_1=\ldots=I_r=I$.
It remains to show that the vertices of the convex hull of $T'$
lie in $H$. For the vertices on $F_l$ this question reduces to the
analogous problem for $GL_l({\bf C})$ and $GL_{e-l}({\bf C})$.
Thus by induction in $e$ it is sufficient to consider the vertex
$({\vert J\vert}k{\vert I\vert}/e,\ldots,
{\vert J\vert}k{\vert I\vert}/e)$. In other words, does 
$V(I)^{\otimes ke}$ contain a one-dimensional representation of
$GL_e({\bf C})$, when $k$ is sufficiently large? 
By explicit application of the Littlewood-Richardson rule 
one can show that already $k=1$ is sufficient.
It is more convenient, however, to consider large $k$ and 
the standard character formula for the multiplicity of the trivial
representation of $SU(e)$ in $V(I)^{\otimes ke}$. When $\chi$ is the
character of the $SU(e)$ representation $V(I)$, this multiplicity is
given by the integral of $\chi^{ke}$ over the Haar measure of $SU(e)$.
With increasing $k$ the multiplicity becomes arbitrarily large, since
the dominant contributions come from small neighborhoods of the
central elements of $SU(e)$.
$\hfill{\Box}$

\bigskip
{\bf Corollary 2.5 :} {\it In the notation of theorem 2.2, ampleness
of $S_{I(l,m,0)}E$ and $\Lambda^m E$ are equivalent. Ampleness
of $S_{I(l,m,s)}E$ implies ampleness of $\Gamma ^\alpha_k E$.}
\bigskip

Before we come to the proof of theorem 2.2, let us
prove a proposition which allows us to get rid of any lower
bound on $k$ when $m=1$.
The procedure may be useful in more general situations,
so we will be somewhat more general than necessary in our context.
\medskip

Let us recall that a line bundle $L$ on a variety $X$ is nef if
$(C.L)\geq 0$ for every curve  $C$ in $X$
and a vector bundle $F$ is nef, when the line
bundle  ${\cal O}_{P(F)}(1)$ is nef on $\bP(F).$  If $F$ is nef then
$S_IF$ is nef for any partition $I$.
\bigskip

{\bf Bloch-Gieseker Lemma [B-G, Lemma.2.1]} {\it  Let $L$ be a line
bundle on
a projective variety $X$
and $d$ a positive integer. Then there exist a projective variety $Y,$ 
a finite
surjective morphism $f:Y\to X$, and a line bundle $M$ on $Y$ such that
 $f^*L\simeq M^d.$ }

\bigskip

{\bf Lemma 2.6 :} {\it Let $X$ be a projective variety. 
Let $E$ and  $F_i$ with $i=1,\ldots,m$ be
vector bundles on $X$ of ranks $e,f_i$
with $E$ ample and the $F_i$ nef. Let $I,J_i$ be non-vanishing
partitions. Then there exist a projective variety $Y,$ a finite
surjective morphism $f:Y\to X$, and ample vector
bundles $E',F'_i$ on $Y$ of ranks $e,f_i$ such that
$f^*(S_IE\otimes_{i=1}^m\,S_{J_i}F_i)=S_I E' \otimes S_{J_i}F'_i$.}
\bigskip

{\it Proof:} If $E$ is ample, $S^kE\otimes (det E)^*$ is ample for
$k>>0$.
We fix $k$ accordingly and use the B-G lemma such that $f^*_{\alpha }
(detE)=M^{k\alpha }$. Then  $E'_{\alpha }=f^*_{\alpha }E\otimes
(M^*)^{\alpha }$
is ample since $S^kE'_{\alpha }$ is.
For $\alpha =\sum_i\vert J_i\vert$ we have
$f^*_{\alpha }(S_IE\otimes_{i=1}^m\, S_{J_i}F_i)=S_IE'_{\alpha }
\otimes S_{J_i}F'_i$ where $F'_i=M^{\vert I\vert}\otimes f^*_{\alpha
}F_i$.
According to Fujita's lemma [Fu], the $F'_i$ are ample, too.
$\hfill{\Box}$

\bigskip

{\bf Corollary 2.7 :} {\it  To prove a vanishing theorem
$H^{p,q}(X,S_IE\otimes_{i=1}^m\, S_{J_i}F_i)=0$
for fixed $p,q,n,e,f_i$
under the assumption  $E$ is ample and the $F_i$ are nef, it
is sufficient to treat the case when
$E$ and the $F_i$ are all ample.}
\bigskip

Indeed, the vanishing of $H^{p,q}(Y,f^*{\cal F})$ implies the vanishing
of $H^{p,q}(X,{\cal F})$
for any  vector bundle ${\cal F}$ on $X$ and any finite surjective
morphism
 $f:Y\to X.$
$\hfill{\Box}$

\bigskip
Now we can show how to get rid of lower bounds on $k$ when $m=1$.
\medskip

{\bf Proposition 2.8 :} {\it Fix $n,p,q,\alpha\in\N$ and $t\in\Z$.
Assume that $H^{p,q}(X,\Gamma^\alpha_k E)$ vanishes for
some $k$, all compact varieties $X$ of dimension $n$ and all ample
vector
bundles $E$ of rank $e=k+t$. Let $k'<k$. Then
$H^{p,q}(X,\Gamma^\alpha_{k'}E')$ vanishes for all ample vector bundles
$E'$ of rank $e'=k'+t$.}
\bigskip

{\it Proof:} For given $E'$, put $E=E'\oplus L^{\oplus(k-k')}$, where
$L$
is any ample line bundle. Since $\Gamma ^\alpha_{k'}E'\otimes L^{k-k'}$
is
a direct summand of  $\Gamma ^\alpha_k E$, we have
$$H^{p,q}(X,\Gamma^\alpha_{k'} E'\otimes L^{k-k'})=0$$
for ample vector bundles $E'$ of rank $e'$ and ample line bundles $L$.
By Corollary 2.7,
this vanishing result remains true, when $L$ is
replaced by the trivial line bundle.
$\hfill{\Box}$

\section{Proof of the theorem 2.2.}

To prepare the proof, we need a lemma and the definition of the Borel-Le
Potier
spectral sequence, which has been made a standard tool
in the derivation of vanishing
theorems [Dem]. \\

Let $E$ be a vector bundle over a compact complex
variety $X$.  Considering subspaces of codimension $s_1,s_2,\ldots$
in the fibres, one obtains the corresponding incomplete flag manifold
$Y=M_{s_1,s_2,\ldots}(E)$, with a natural projection  $\pi:Y\rightarrow
X$.
Let $Q_k$ with $rank(Q_1)=s_1$, $rank(Q_2)=s_2-s_1, \ldots$ be the
corresponding
canonical quotient bundles over $Y$.
Define the line bundles
$$Q^{i_1,i_2,\ldots}=\bigotimes_k\ det(Q_k)^{i_k}\ .$$

{\bf Lemma 3.1 :} {\it Let $I$ be a partition of the form
$$I=(\underbrace {i_1,\ldots ,i_1}_{s_1\mit\ times},
\underbrace {i_2, \ldots ,i_2}_{(s_2-s_1)\mit\ times},\ldots)$$
such that $S_IE$ is ample over $X$. Then $Q^{i_1,i_2,\ldots}$ is ample
over $Y$.}
\medskip

{\it Proof:} See [Dem], Lemmata 2.11, 4.1.
$\hfill{\Box}$

\bigskip

The projection $\pi$ yields a filtration of the bundle
$\Omega^P_Y$ of exterior differential forms of degree $P$, namely

$$F^p(\Omega^P_Y)=\pi^*\Omega^P_X \wedge \Omega^{P-p}_Y.$$
The corresponding graded bundle is given by
$$F^p(\Omega^P_Y)/F^{p+1}(\Omega^P_Y)=
\pi^*\Omega^P_X\otimes \Omega^{P-p}_{Y/X},$$
where $\Omega^{P-p}_{Y/X}$ is the bundle of relative
differential forms of degree $P-p$.

For a given line bundle ${\cal L}$ over $Y$
the corresponding filtration of $\Omega^P_Y\otimes {\cal L}$
yields a spectral sequence which abuts on $H^{P,q}(Y,{\cal L})$.
It has been named
Borel-Le Potier  spectral sequence by Demailly.
It is given by the data $X,Y,{\cal L},P$ and will be denoted by $^PE$.
Its $E_1$-term

$$^PE_1^{p,q-p}=
H^q(Y,\pi^*(\Omega^p_X)\otimes\Omega^{P-p}_{Y/X}\otimes {\cal L})$$
can be calculated  via the Leray spectral sequence.
\bigskip

Now let us come to the proof.
For convenience, we add an index $0$ to the variables $l,s,\alpha,
p,q,r$ in theorem 2.2, such that without index
they can be used as free variables in the proof.
Thus we  prove a vanishing theorem
for $H^{p_0,q_0}(X,\Gamma_k^{\alpha_0}E)$ where
$0\leq p_0,q_0\leq n$, $0\leq \alpha_0<k$
and we define
$$r_0=\delta \left(n-p_0+{m\choose 2}\right)
\ .$$

We will work by induction in the finite set
$$B=\{(\alpha,p,q)\in \N^3\,\vert\, k-e\leq \alpha\leq \alpha_0+r_0-m,\,
p_0+Q(\alpha)\leq p\leq n,\, q_0+Q(\alpha)\leq q\leq n\},$$
where
$$Q(\alpha)=r_0(\alpha-\alpha_0)-{|\alpha-\alpha_0|\choose 2}\ .$$
By a suitable ordering of $B$, we shall prove inductively that
$H^{p,q}(X,\Gamma^\alpha_k E)=0$ for all $(\alpha,p,q)\in B$, in
particular for its maximal element $(\alpha_0,p_0,q_0)$.
\bigskip

For any $(\alpha,p,q)\in B$ consider the incomplete flag variety
$Y=M_{s,r}(E)$ with $r=r_0+\alpha_0-\alpha$,\ $l=[k/r]$ and $s=k-lr$.
This variety reduces to the Grassmannian $G_r(E)$ when $s=0$.
For the line bundle ${\cal L}= Q^{l+1,l}$ and
$$P=p-r\alpha+(l-1){r+1\choose 2}+rs-{s\choose 2}$$
consider the corresponding Borel-Le Potier spectral sequence $^PE$.
By the Leray spectral sequence,
$$^PE_1^{p',q'-p'}=\bigoplus_{j\in\N}
 H^{p',q'-j}(X,R^j\pi_*(\Omega_{Y/X}^{P-p'}\otimes {\cal L}))\ .$$

We shall see that the r.h.s. can be evaluated in terms of hook
Schur functors, if either $r=1$ or  $k\geq n-p+r\alpha+l$.
For $m=1$ we can assume that this lower bound is satisfied for all
elements of $B$, since proposition 2.8 allows to get rid of it
subsequently. Otherwise we have to show that
$k>k(\alpha)$, where
$$k(\alpha)=(A(\alpha)-1)(1+(r_0+\alpha_0-\alpha-1)^{-1})$$
and
$$A(\alpha)=n-p_0-Q(\alpha)+r\alpha\ .$$
For $r_0+\alpha_0-\alpha=1$ we put $k(\alpha)=-\infty$.

When $k(\alpha)$ is extended to real argument, it becomes
a decreasing function of $\alpha$ in the interval
$\alpha_0\leq \alpha\leq r_0+\alpha_0-1$,
$1\leq \alpha$, since
$$n-p_0+r_0\alpha_0+{\alpha_0+1\choose 2}<{r_0+\alpha_0+1\choose 2}\ .$$
For $\alpha_0=0$, one checks easily that $k(0)\geq k(1)$, such that
$k>k(0)$ implies $k>k(\alpha)$ for all $\alpha$ occuring in $B$.

Thus $k>k(0)$ implies $k>k(\alpha)$ for all $\alpha$ occuring in $B$,
if $\alpha_0=0$.
For $1\leq \alpha\leq \alpha_0$ one has $k(1)\geq k(0)$. Moreover,
the function $k(\alpha)$ is
non-increasing for $\alpha\geq 2$, $r_0+\alpha_0\geq 6$ and for
$\alpha\geq 1$, $r_0+\alpha_0\geq 7$.
Thus the function $[k(\alpha)]$, $\alpha$ occuring in $B$,
takes its maximum at $\alpha=1$ or $\alpha=k-e$, except for
a small number
of cases for which the maximum of $[k(\alpha)]$ turns out to be
$[k(1)]+1$. For $m>1$ this only happens for
$m=2$, $n-p_0=1$, $\alpha_0=2$.

Let
$$\chi=(l-1){r\choose 2}+rs-{s+1\choose 2}-(r-1)\alpha\ .$$
Since either $r=1$ or  $k\geq n-p+r\alpha+l$
one has according to [Ma], Prop. 3 and  Lemma 3,
$$^PE_1^{p',q'+\chi-p'}=\bigoplus_{\beta=0}^{[(k-1)/r]}
n_s(\sigma+p-p'+r(\beta-\alpha))
\,H^{p',q'+p'-p+\alpha-\beta}(X,\Gamma_k^\beta E),$$
where $\sigma=rs-{s\choose 2}$ and the multiplicity function $n_s$ is
generated by
$$\sum_{a,s} n_s(a)x^az^s=\prod_{i=1}^r(1+x^{r+1-i}z)\ .$$
Note that $n_s(a)=0$ for $a>\sigma$ and $n_s(\sigma)=1$.

Manivel states his result under the additional condition $e\geq k$,
which is
not necessary and not used in the proof.
For $k>e$ the lower bound of the $\beta$ summation can be replaced
by $k-e$. The upper bound is given as $l$ in Manivel's statement.
This is correct for $s>0$ but has to be reduced by 1 for $s=0$, as
is clear from his proof. A priori, one only expects the occurence
of Schur functors $\Gamma_k^\beta$ which are dominated
by $S_{I(l,r,s)}$, which indeed is true for the correct upper bound.
\bigskip
Note that $^PE_1^{p,q+\chi-p}$ contains a direct summand
$H^{p,q}(X,\Gamma ^\alpha_k E)$.
\medskip

In the Borel-Le Potier spectral sequence consider any non-vanishing
morphism
$$^PE_M^{p,q+\chi-p}\longrightarrow\ ^PE_M^{p+M,q+\chi+1-p-M},$$
with $M>0$, or
$$^PE_{-M}^{p+M,q+\chi-1-p-M}\longrightarrow\ ^PE_{-M}^{p,q+\chi-p},$$
with $M<0$. Then $^PE_1^{p+M,q+\chi+sgn(M)-p-M}$ is a direct sum of
terms
$H^{p',q'}(X,\Gamma_k^\beta E)$ with $0\leq p',q'\leq n$ such that
$p'-p\geq r(\beta-\alpha)$ and
$$q'-q=sgn(M)+\alpha-\beta+M\geq (r-1)(\beta-\alpha)+sgn(\beta-\alpha)\
.$$
\smallskip

From these conditions one obtains $(\beta,p',q')\in B$. Indeed, the
conditions
for $p',q'$ follow from the easily verified relations
$Q(\beta)-Q(\alpha)\leq r(\beta-\alpha)$ and
$Q(\beta)-Q(\alpha)\leq (r-1)(\beta-\alpha)+sgn(\beta-\alpha)$.
We still have to check that $\beta\leq \alpha_0+r_0-m.$
Note that we already have shown
$p'\geq Q(\beta)+p_0$. By the definition of $\delta$ and
$Q$ one finds
$Q(r_0+\alpha_0-m+1)+p_0>n$. Thus
for $\beta= \alpha_0+r_0-m+1$
and {\it a fortiori} for $\beta> \alpha_0+r_0-m+1$ one has
$p'>n$, which is excluded.

\bigskip

To prove the desired vanishing theorem by induction, we still need
$(\beta,p',q')<(\alpha,p,q)$ in a suitable ordering. We use
the variable
$$L(\alpha,p)=2(n-p)+\alpha(2(r_0+\alpha_0)-\alpha)\ .$$
A short calculation yields $L(\beta,p')\leq L(\alpha,p)-1$ for any
of the morphisms under consideration.
\bigskip

Assume that $H^{p',q'}(X,\Gamma_k^\beta E)=0$ for all 
$(\beta,p',q')\in B$
with $L(\beta,p')< L(\alpha,p)$. Then
$$^PE_1^{p,q+\chi-p}=\ ^PE_\infty^{p,q+\chi-p}\ ,$$
since the morphisms adjacent to $^PE_M^{p,q+\chi-p}$
all vanish.
By the Kodaira-Akizuki-Nakano vanishing theorem the r.h.s. vanishes
when $Q^{l+1,l}$ is ample and $P+q+\chi>dim(Y)$.
By lemma 2.4 and lemma 3.1, the ampleness of $S_{I(l,m,s)}E$
guarantees indeed that $Q^{l+1,l}$ is ample for all $(\alpha,p,q)\in B$.

Since $dim(Y)=n+r(e-r)+s(r-s)$
and $p+q-n\geq p_0+q_0-n+2Q(\alpha)$ the condition  $P+q+\chi>dim(Y)$
reduces to
$$p_0+q_0-n+\alpha(e-k+\alpha)+|\alpha-\alpha_0|+\alpha-\alpha_0
>(r_0+\alpha_0)(e-k+\alpha_0)+\alpha_0(r_0-1)\ .$$
By assumption, $p_0+q_0-n$ is strictly greater than the r.h.s.
Since $\alpha(e-k+\alpha)$ and $|\alpha-\alpha_0|+\alpha-\alpha_0$ are
both
non-negative (they only vanish for the minimal value of $\alpha$
occuring in $B$), the conditions of the Kodaira-Akizuki-Nakano
vanishing theorem are satisfied. Thus $^PE_1^{p,q+\chi-p}$ and its
direct summand $H^{p,q}(X,\Gamma ^\alpha_k E)$ do vanish. In
particular, this is true for
$(\alpha_0,p_0,q_0)\in B$, as was to be shown.
\medskip

When we put
$$r_0=\delta \left(n-q_0+{m\choose 2}\right)\ ,$$
one has to replace $p$ by $q$
in two passages of the proof.
The first one occurs when one shows that $(\beta,p',q')\in B$.
One finds $Q(r_0+\alpha_0-m+1)+q_0>n$.
One now concludes for $\beta\geq \alpha_0+r_0-m+1$
that $q'>n$, which is excluded.

The second replacement concerns the induction variable
which now becomes $L(\alpha,q)$.
For any of the morphisms under consideration one finds
$L(\beta,q')\leq L(\alpha,q)-1$.
Apart from these single letter replacements, the proof goes through
verbally.
$\hfill{\Box}$
\medskip

Remarks: The condition $k\geq k_0$ is necessary to exclude the 
occurence of
non-hook partitions in the spectral sequence. It should be possible
to weaken the lower bound by suitable vanishing theorems for these 
extra terms.
\bigskip

The ampleness condition of theorem 2.2 is not the only possible one.
When
$E$ allows the extraction of a sufficiently positive line bundle, one
can
use the idea of proposition 2.8 to increase $k$ to a value divisible
by $m$.
\bigskip

{\bf Acknowlegements:} Part of the work was done at the ICTP Trieste.
\bigskip

{\bf References.}

\begin{itemize}

\item{[B-G]} S. Bloch, D. Gieseker, {\it The positivity of the Chern
classes of an ample vector bundle,} Invent. Math. {\bf 12} (1971)
112-117.

\item{[Bu]} A.S. Buch, {\it The Saturation Conjecture (after A. Knutson
and T. Tao)}, preprint, November 1998.

\item{[Dem]} J.P. Demailly,   {\it Vanishing theorems for powers
 of an ample vector bundle,} Invent. Math. {\bf 91} (1988), 203-220.

\item{[Fu]} T. Fujita,  {\it Semipositive Line Bundles}, J. Fac.
Sci. Univ. Tokyo Sect. IA
Math. {\bf 30} (1983), 353-378.

\item{[Gr]} P.A. Griffiths,  {\it Hermitian diffferential geometry, 
Chern classes and positive vector bundles}, Global analysis,
paper in honor of Kodaira (1969), 185-251.

\item{[K-T] A. Knutson and T. Tao {\it The honeycomb model of the
Berenstein-Zelevinsky polytope I: Klyachko's saturation conjecture},
preprint, 1998.

\item{[La]}  F. Laytimi, {\it On Degeneracy Loci}, Int. J. Math. {\bf 7}
(1996), 745-754.

\item{[LeP]} J. Le Potier,   {\it Annulation de la cohomology a valeurs
dans un fibre vectoriel holomorphe positif de rang quelconque},
Math. Ann. {\bf 218} (1975), 35-53.

\item{[Ma]} L. Manivel,  {\it Un th\'eor\`eme d'annulation pour les
puissances ext\'erieures d'un fibr\'e ample}, J. reine angew. Math. 
{\bf 422}
(1991), 91-116.

\item{[Na]} W. Nahm,  {\it A vanishing theorem for products of exterior
and
symmetric powers},
Preprint. (1995), Univ. Bonn. Germany.

\item{[Ze]} A. Zelevinsky, {\it Littlewood-Richardson semigroups},
MSRI Preprint, 1997-044.
}
\end{itemize}

\end{document}